\documentclass[10pt,a4paper]{article}

\usepackage{amsmath}  
\usepackage{amssymb}  
\usepackage{amsthm}   
\usepackage{amsfonts}
\usepackage{mathrsfs}  

\newcommand{\N}{\mathbb{N}}

\newcommand{\Hom}[0]{\operatorname{Hom}}

\newcommand{\Supp}[0]{\operatorname{Supp}}
\newcommand{\Ass}[0]{\operatorname{Ass}}
\newcommand{\Spec}[0]{\operatorname{Spec}}

\newcommand{\dirlim}[0]{\operatorname{dirlim}}
\newcommand{\Ext}[0]{\operatorname{Ext}}
\newcommand{\height}[0]{\operatorname{height}}
\newcommand{\injdim}[0]{\operatorname{injdim}}
\newcommand{\depth}[0]{\operatorname{depth}}
\newcommand{\rank}[0]{\operatorname{rank}}
\newcommand{\Ann}[0]{\operatorname{Ann}}
\newcommand{\cd}[0]{\operatorname{cd}}
\newcommand{\bigheight}[0]{\operatorname{beight}}

\newtheorem{satz}{Theorem}[section]
\newtheorem{lemma}[satz]{Lemma}
\newtheorem{corollary}[satz]{Corollary}
\newtheorem{remark}[satz]{Remark}
\newtheorem{question}[satz]{Question}
\newtheorem{example}[satz]{Example}
\newtheorem{examples}[satz]{Examples}

\newtheorem{summary}[satz]{Summary}
\newtheorem*{satzo}{Theorem}
\newtheorem*{questiono}{Question}

\theoremstyle{definition}

\theoremstyle{plain}

\theoremstyle{remark}

\author{M. Hellus}

\title{A Note on the injective dimension of local cohomology modules}
\date{\today}   

\begin{document}

\maketitle

\begin{abstract}
For a noetherian ring $R$ we call an $R$-module $M$ cofinite if
there exists an ideal $I$ of $R$ such that $M$ is $I$-cofinite; we
show that every cofinite module $M$ satisfies $\dim _R(M)\leq
\injdim _R(M)$. As an application we study the question which local
cohomology modules $H^i_I(R)$ satisfy $\injdim _R(H^i_I(R))=\dim
_R(H^i_I(R))$. There are two situations where the answer is
positive. On the other hand we present two counter-examples, the
failure in these two examples coming from different reasons.
\end{abstract}

\section{Introduction}

Let $I$ be an ideal of a noetherian ring $R$. By definition, an
$R$-module $M$ is $I$-cofinite if $\Supp _R(M)\subseteq {\cal V}(I)$
and if all $\Ext ^i_R(R/I,M)$ are finitely generated. Note that, for
local $(R,m)$, $m$-cofiniteness is the same as artinianness. The
concept of $I$-cofiniteness was introduced by Hartshorne
(\cite{hartshorne70}); see
\cite{delfino97,huneke91,kawasaki98,marley02,yoshida97} for more
material on cofiniteness.

Let $M$ be a module over a ring $R$. An injective resolution of $M$
is an exact sequence \[ 0\to M\to E^0\to E^1\to \dots \ \ ,\] where
all $E^i$ are injective $R$-modules. By definition, $\injdim
_R(M)\leq n$, if there is an injective resolution of $M$ such that
$0=E^{n+1}=E^{n+2}=\dots $.

Every finitely generated $M$ over a noetherian ring $R$ satisfies
$\dim _R(M)\buildrel (1)\over \leq \injdim _R(M)$; if $R$ is local
and $\injdim _R(M)<\infty $ holds, one has $\injdim _R(M)\buildrel
(2)\over =\depth (R)$. Both $(1)$ and $(2)$ are known as Bass
formula. In this work (theorem \ref{genBass}) we generalize $(1)$ to
the class of all $R$-modules which are $I$-cofinite for at least one
ideal $I$ of $R$ (we call these modules cofinite modules), where $R$
is a noetherian ring. We show (remark \ref{fails}) that $(2)$ does
not generalize to the class of all cofinite modules and explain why
this generalization fails.

In the sequel we prove some consequences of theorem \ref{genBass}
for local cohomology modules; in particular we are dealing with the
following

\begin{questiono}
When does
\[\injdim _R(H^i_I(R))=\dim _R(H^i_I(R))\] hold?
\end{questiono}

(Here $I$ is an ideal of a noetherian ring $R$ and $H^i_I(R)$ is the
$i$-th local cohomology of $R$ supported in $I$.) There is a
positive answer in the following two situations (see Corollary
\ref{combination} and remark \ref{Gor}):
\begin{itemize}
  \item $R$ is a noetherian local regular ring containing a field and $H^i_I(R)$ is $I$-cofinite.
  \item $R$ is a noetherian local Gorenstein ring $R$ and $H^n_I(R)=0$ for every $n\neq i$.
\end{itemize}
Finally we prove that the answer to the above question is negative
in general. More precisely, we present two counter-examples which
are different in nature (example \ref{mixed} and the last example
from \ref{series}). In \ref{summary} we summarize answers to the
above question.

By $E_R(M)$ we denote a fixed $R$-injective hull of an $R$-module
$M$.

\section{Results}

\subsection{A positive result}

The following statement is contained in \cite[Corollary
1]{delfino97}, for the sake of completeness we present a simple
proof:

\begin{lemma}
\label{onExt}
Let $R$ be a noetherian ring, $I\subseteq R$ an ideal
and $M$ an $I$-cofinite $R$-module. Then, for every ideal $J$ of $R$
containing $I$ and for every $l\in \N $, the module $\Ext
^l_R(R/J,M)$ is finitely generated.
\end{lemma}

{\it Proof. }By induction on $l$: The case $l=0$ is trivial, because
for every ideal $J$ of $R$ containing $I$ there is an inclusion
$\Hom _R(R/J,M)\subseteq \Hom _R(R/I,M)$ and the latter module is
finitely generated. We assume $l>0$ and that the statement is true
for smaller $l$: We choose $x_1,\dots ,x_n\in R$ such that
$J=I+(x_1,\dots ,x_n)R$ and prove our induction hypothesis by
induction on $n$: In the case $n=0$ there is nothing to prove, we
assume $n>0$ and that the statement holds for smaller $n$: Consider
the short exact sequence \[ 0\to R/(I_{n-1}:_Rx_n)\buildrel x_n
\over \to R/I_{n-1}\to R/J\to 0\ \ ,\] (we define $I_s=I+(x_1,\dots
,x_s)R$ for every $s\in \{ 1,\dots ,n\} $) where $x_n$ denotes
multiplication by $x_n$ on $R/I_{n-1}$. We get an exact sequence
\[ \Ext ^{l-1}_R(R/(I_{n-1}:_Rx_n),M)\to \Ext ^l_R(R/J,M)\to \Ext ^l_R(R/I_{n-1},M)\ \ .\] The latter module in this sequence
is finitely generated by induction on $n$ and the first module is
finitely generated by induction on $l$. Therefore, the second module
is finitely generated; both inductions are completed.\hfill $\square
$

\begin{remark}
\label{fails} If $(R,m)$ is a noetherian local ring and $M$ is a
finitely generated $R$-module of finite injective dimension then
$\injdim _R(M)\buildrel (2)\over =\depth (R)$. This is not true for
every cofinite $R$-module $M$, here is a counterexample: Let
$R=k[[X]]$ be a formal power series ring in one variable over a
field $k$; set $M=k[X^{-1}]=k\cdot 1\oplus k\cdot X^{-1}\oplus \dots
$. It is well-known that $M$, together with its natural $R$-module
structure, is an $R$-injective hull of $k$. Clearly, $M$ is
artinian, i. e. $XR$-cofinite and we have \[ \injdim _R(M)=0\neq
1=\depth (R)\ \ .\] More generally, remark \ref{Gor} below will
present a class of local cohomology modules $M$ for which \[ \injdim
_R(M)=\dim _R(M)\] holds, and for which $\dim _R(M)$ is not equal to
$\depth (R)$ in general.
\end{remark}

If $R$ is a noetherian ring, every finitely generated $R$-module $M$
satisfies \[ \dim _R(M)\leq \injdim_R(M)\ \ .\] The same formula
holds for any cofinite $R$-module $M$ (by definition, $M$ is
cofinite if there exists an ideal $I$ of $R$ such that $M$ is
$I$-cofinite):

\begin{satz}
\label{genBass}Let $R$ be a noetherian ring and $M$ a cofinite
$R$-module. Then
\[ \dim _R(M)\leq \injdim _R(M)\] holds.
\end{satz}
{\it Proof. }Let $I$ be an ideal of $R$ such that $M$ is
$I$-cofinite. It suffices to show the following two statements (by
$\mu _k(p,M):=\rank _{R_p/pR_p}(\Ext ^k_{R_p}(R_p/pR_p,M_p))$ we
denote the $k$-th Bass number of $M$ with respect to a prime ideal
$p$ of $R$):
\par
(a) If $p$ is minimal in $\Supp _R(M)$ then $\mu _0(p,M)\neq 0$.
\par
(b) If $p\subsetneq p^\prime $ are prime ideals of $R$ with no prime
in between such that $\mu _k(p,M)\neq 0$ then $\mu _{k+1}(p^\prime
,M)\neq 0$:
\par
Let $0\to M\to E^\bullet $ be a minimal injective resolution of $M$.
Statement (a) is clear, because for every $p$ minimal in $\Supp
_R(M)$ one has $p\in \Ass _R(M)=\Ass _R(E^0)$, i. e. $\mu
_0(p,M)\neq 0$. Proof of (b): Localizing at $p^\prime $ allows us to
assume that $(R,m)$ is local and $p^\prime =m$. Choose any $x\in
m\setminus p$. The short exact sequence \[ 0\to R/p\buildrel x\over
\to R/p\to R/(p+xR)\to 0\] induces an exact sequence \[ \Ext
^k_R(R/p,M)\buildrel x\over \to \Ext ^k_R(R/p,M)\to \Ext
^{k+1}_R(R/(p+xR),M)\ \ .\] Now, $\mu _k(p,M)\neq 0$ implies $p\in
\Supp _R(M)$ (note that one has $\Supp _R(E^l)\subseteq \Supp _R(M)$
for every $l$ as $E^\bullet $ is a minimal injective resolution of
$M$) and, therefore, $p\supseteq I$; by lemma \ref{onExt}, we
conclude that $\Ext ^k_R(R/p,M)(\neq 0)$ is a finitely generated
$R$-module. This fact, together with our last exact sequence and the
lemma of Nakayama implies that $\Ext ^{k+1}_R(R/(p+xR),M)\neq 0$.
The non-zero $R$-module $R/(p+xR)$ has finite length, and therefore
it is easy to see that $\Ext ^{k+1}_R(R/m,M)\neq 0$, i. e. $\mu
_{k+1}(m,M)\neq 0$.\hfill $\square $

\begin{corollary}
\label{c1c2} Let $R$ be a noetherian ring, $I$ an ideal of $R$ and
$M$ a finitely generated $R$-module. Then $H^l_I(M)$ is $I$-cofinite
for every $l\in \N$ if one of the following two conditions holds:

$(c_1)$ $\dim (R/I)=1$.

$(c_2)$ There exists $p\in \N$ such that $H^i_I(M)$ is finitely
generated for every $i\neq p$ (see remark \ref{c2} on this
condition; also note that it would suffice to assume that $H^i_I(M)$
is $I$-cofinite for all $i\neq p$).

Consequently, by theorem \ref{genBass}, one has
\[ \dim_R(H^l_I(M))\leq \injdim _R(H^l_I(M))\]
for every natural $l$.
\end{corollary}

{\it Proof. }This follows immediately from theorem \ref{genBass} and
the fact that $H^l_I(M)$ is $I$-cofinite for every $l\in \N$ if
$(c_1)$ or $(c_2)$ holds; this was shown in \cite[Theorem
1]{delfino97} resp. \cite[Theorem 1.1]{yoshida97} in the case of
$(c_1)$ and in \cite[Proposition 2.5]{marley02} in the case of
$(c_2)$.\hfill $\square $

\begin{remark}
\label{c2} In the situation of corollary \ref{c1c2} set \[
h:=\height ((I+\Ann_R(M))/\Ann _R(M))\] and \[ c:=\cd (I,M)=\cd
(I,R/\Ann _R(M))\ \ .\] Note that $\cd (I,M):=\inf \{ l\in \N \vert
 0=H^{l+1}_I(M)=H^{l+2}_I(M)=\dots \} $ is the cohomological
dimension of $I$ on $M$; also note that $0\leq h\leq c$. One has
\[ c>0\Rightarrow H^c_I(M)\hbox { is not finitely generated}\]
(proof: By localizing in any prime ideal of the support of
$H^c_I(M)$ we may assume that $(R,m)$ is local; the functor $H^c_I$
is right exact on modules whose support is contained in
$V(Ann_R(M))$ and therefore one has $H^c_I(M)\otimes
_RR/m=H^c_I(M/mM)=0$, as $c$ is positive; by Nakayama, $H^c_I(M)$ is
not finite) and
\[ h>0\Rightarrow H^h_I(M)\hbox { is not finitely generated}\]
(proof: We assume that $H^h_I(M)$ is finitely generated. Let $p$ be
a prime ideal of $R$ containing $I+\Ann _R(M)$ such that $\height
(p/\Ann _R(M))=h$. Then $IR_p/\Ann_R(M)R_p$ is
$pR_p/\Ann_R(M)R_p$-primary and we conclude that
\[ H^h_I(M)_p=H^h_{pR_p}(M_p)\]
is finite as an $R_p$-module; by construction, $\dim _{R_p}(M_p)=h$
and our above result on $c$ shows that $H^h_{pR_p}(M_p)$ is not
finite, contradiction). Thus, if $h>0$, the equality $p=h=c$ is
necessary for condition $(c_2)$ of corollary \ref{c1c2}.
\end{remark}
\subsection{An application and some examples}

\begin{corollary}
\label{combination} Let $R$ be a noetherian local regular ring
containing a field and $I$ an ideal of $R$ such that $H^l_I(R)$ is
$I$-cofinite for some natural number $l$. Then
\[ \injdim _R(H^l_I(R))=\dim _R(H^l_I(R))\]
holds. In particular, one has $\injdim _R(H^l_I(R))=\dim
_R(H^l_I(R))$ if one of the following conditions is satisfied:

(i) $\dim (R/I)\leq 1$.

(ii) The residue field $R/m$ (where $m$ denotes the maximal ideal of
$R$) is separably closed, $\dim(R/I)=2$, $I$ is equidimensional and
$\Spec(R/I)\setminus\{ m\} $ is connected.

(iii) $I$ is a set-theoretic complete intersection.

(iv) The characteristic of $R$ and $R/I$ is Cohen-Macaulay.
\end{corollary}

{\it Proof. }The first statement follows immediately from theorem
\ref{genBass} together with the following theorem, which is due to
Lyubeznik (\cite{lyubeznik93,lyubeznik97}):

\begin{satzo}
Let $R$ be a regular ring containing a field and $I$ an ideal of
$R$. Then, for every $l\in \N $, \[ \injdim _R(H^l_I(R))\leq \dim
_R(H^l_I(R))\] holds.
\end{satzo}

In the second part it is true in all four cases that there is only
one non-vanishing local cohomology module of $R$ with support in $I$
(which is then necessarily $I$-cofinite by \cite [Prop.
2.5]{marley02}); this is clear for the cases (i) and (iii). In the
case of (ii), it follows from \cite[Theorem 2.9]{huneke90} and in
the case of (iv) from the identity
\[ H^l_I(R)=\dirlim _{e\in \N}\Ext ^l_R(R/I^{[p^e]},R)\]
(where $I^{[p^e]}$ denotes the $p^e$-th Frobenius power of $I$ and
$l$ is an arbitrary natural number) together with the flatness of
the Frobenius map over regular local rings of positive
characteristic.

\hfill $\square $

Under suitable assumptions, the result of corollary
\ref{combination} holds for more general rings:

\begin{remark}
\label{Gor}
Let $I$ be an ideal of a Gorenstein ring $R$ and $l\in \N
$ such that \[ H^i_I(R)=0\ \ (i\neq l)\] holds. Then the formula
\[ \injdim _R(H^l_I(R))=\dim _R(H^l_I(R))\ (=\dim (R/I))\] holds.
\end{remark}
{\it Proof. }Follows e.~g. from \cite[Prop. 3.1]{huneke91}.\hfill
$\square $

\begin{question}
\label{when} When does $\injdim _R(H_I^i(R))=\dim H_I^i(R)$ hold? In
the situations of corollary \ref{combination} and remark \ref{Gor}
this question has a positive answer.
\end{question}

But corollary \ref{combination} in general becomes false if one does
not assume that $H^l_I(R)$ is $I$-cofinite, as the following example
shows:

\begin{example}
\label{mixed} Let $R=k[[x,y,z]]$ be a power series ring over a field
$k$ and let $I$ be the ideal $(xy,xz)R$ of $R$. Because of $I=xR\cap
(y,z)R$ we have the Mayer-Vietoris sequence
\[ 0\to H^2_{(y,z)R}(R)\to H:=H^2_I(R)\to E\to 0\]
($E$ is a fixed $R$-injective hull of $k$) and hence an exact
sequence
\[ \Hom _R(R/(y,z)R,H^2_{(y,z)R}(R))\to \Hom _R(R/(y,z)R,H)\to \] \[ \to \Hom
_R(R/(y,z)R,E)\to \Ext ^1_R(R/(y,z)R,H^2_{(y,z)R}(R))\ .\] The first
and the last term of the previous sequence are both finitely
generated, since $H^2_{(y,z)R}(R)$ is $(y,z)R$-cofinite (the last
term is zero in fact). Therefore, $\Hom _R(R/(y,z)R,H)$ is finite if
and only if $\Hom _R(R/(y,z)R,E)$ is finite; but $\Hom
_R(R/(y,z)R,E)$ is isomorphic to $E_{k[[x]]}(k)$ (with an $R$-module
structure induced by $R\to R/(y,z)R\cong k[[x]]$) and it is well
known that the latter module is not finitely generated. We conclude
that $\Hom _R(R/(y,z)R,H)$ is not finitely generated; in particular
$\Hom _R(R/I,H)$ is not finitely generated, $H$ is not $I$-cofinite.

By using Cech cohomology, it is clear that $H$ is the cokernel of
the natural map
\[ R_{xy}\oplus R_{xz}\to R_{xyz}\ \ .\]
Thus, by using Cech cohomology also for $H^2_{(y,z)R}(R)$,
\[
H=H^2_{(y,z)R}(R)_x=H^2_{(y,z)R_{(y,z)R}}(R_{(y,z)R})=E_R(R/(y,z)R)\]
(the second equality follows because $(y,z)R_x$ is a maximal ideal
of the localized ring $R_x$ and $H^2_{(y,z)R}(R)_x$ is an
$R_{(y,z)R}$-module). $H$ is injective, we have
\[ \injdim _R(H)=0\neq 1=\dim _R(H)\ \ .\]
\end{example}

It seems that the negative answer to question \ref{when} in the
previous example comes from the fact that the primary components of
the ideal $I$ have different dimension; but we will see that the
answer is in general negative also if all components of $I$ have the
same dimension. We start with some general facts:

\begin{remark}
\label{nichtCof} There is a result from Huneke and Koh
(\cite[theorem 2.3(ii)]{huneke91}) resp. Lyubeznik (\cite[corollary
3.5]{lyubeznik93}) which says that if $I$ is an ideal of a regular
ring $R$ containing a field, then, for $l>\bigheight(I)$, $\Hom
_R(R/I,H^l_I(R))$ is finitely generated only if $H^l_I(R)$ is zero.
Thus, in order to find examples of non-$I$-cofinite local cohomology
modules $H^l_I(R)$ one should consider ideals where $\cd (I,R)$ is
big.

Faltings (\cite{faltings80}) proved upper bounds for the
cohomological dimension $\cd (I,R)$ of an ideal $I$ of a regular
local ring $R$ containing a field. Lyubeznik presents examples where
this upper bounds are actually obtained (\cite{lyubeznik85} and
\cite[Corollary 5.3]{huneke90}):
\end{remark}

\begin{satzo} Let $(R,m)$ be an excellent regular $d$-dimensional local ring
containing a field. Let $s$ and $b$ be positive natural numbers such
that $s=\lceil (d-1)/b\rceil $, where $\lceil \alpha \rceil $is the
maximum integer which does not exceed the real number $\alpha $. Let
$I_0,\dots ,I_s$ be ideals of bigheight $b$ such that $I_0+\dots
+I_s$ is $m$-primary. Set $I=\bigcap _{0\leq j\leq s}I_j$. Then $\cd
(I,R)=d-s$ and this is precisely the maximal possible value (i. e.
the general upper bound on $\cd (I,R)$ from \cite{faltings80}).
\end{satzo}

We investigate, with respect to question \ref{when}, non-trivial
examples where $s$ and $b$ are small (note that cases where the
local cohomology module is cofinite are clear by Corollary
\ref{combination}); our last example will show that not all local
cohomology modules $H$ arising from Lyubeznik's result above satisfy
$\injdim _R(H)=\dim _R(H)$: {\sloppy
\begin{examples}
\label{series} Let $k$ be a field and $R=k[[x_1,\dots ,x_d]]$ a
formal power series algebra over $k$ in $d$ variables.
\begin{itemize}
  \item $s=2,b=2,d=5$. We take $I:=(x_1,x_2)R\cap (x_3,x_4)R\cap
  (x_5,x_1)R$ and $H:=H^3_I(R)$. Because of remark \ref{nichtCof}, $H$ is not $I$-cofinite (it is not zero because, by Lyubeznik's result mentioned above, $\cd
  (R,I)=5-s=3$).
  We claim that $H$ has a minimal
  injective resolution of the form
  \[ 0\to H\to E_R(R/(x_1,x_2,x_3,x_4)R)\oplus
  E_R(R/(x_1,x_3,x_4,x_5)R)\to \] \[\to E_R(R/m)\to 0\ :\]
  A part of a Mayer-Vietoris sequence with respect to the ideal
  $(x_1,x_2)R\cap (x_5,x_1)R$ and $(x_3,x_4)R$ is:
  \[ H^3_{(x_1,x_2)R\cap (x_5,x_1)R}(R)\to H\to \]
  \[ \to H^4_{(x_1,x_2,x_3,x_4)R\cap (x_1,x_3,x_4,x_5)R}(R)\to
  H^4_{(x_1,x_2)R\cap (x_5,x_1)R}(R)\]
  By another Mayer-Vietoris sequence argument, we see that in this
  sequence the first and last term vanish, we get
  \[ H=H^4_{(x_1,x_2,x_3,x_4)R\cap (x_1,x_3,x_4,x_5)R}(R)\]
  For later use we note $\{ (x_1,x_2,x_3,x_4)R,(x_1,x_3,x_4,x_5)R\} \buildrel (+)\over \subseteq
  \Ass _R(H)$. The Cech complex of $H$ with respect to $x_2,x_5$ has
  the form
  \[ 0\to H\to H_{x_2}\oplus H_{x_5}\to H_{x_2\cdot x_5}=0\ \ .\]
  The spectral sequence belonging to the composed functors
  \[ \Gamma
  _{(x_2,x_5})\circ \Gamma _{(x_1,x_2,x_3,x_4)R\cap
  (x_1,x_3,x_4,x_5)R}\]
  shows that $\Gamma _{(x_2,x_5)R}(H)=0$ and
  $H^1_{(x_2,x_5)R}(H)=E_R(R/m)$. On the other hand, for
  $p:=(x_1,x_3,x_4,x_5)R$ we have
  \[ H_{x_2}=H^4_{pR_{x_2}}(R_{x_2})=H^4_{pR_p}(R_p)=E_R(R/(x_1,x_3,x_4,x_5)R)\]
  (because $H^4_{pR_{x_2}}(R_{x_2})$ has a natural $R_p$-module structure). Similarly, $H_{x_5}=E_R(R/(x_1,x_2,x_3,x_4))R$.
  Together with $(+)$ our claim follows. Thus,
  \[ \injdim _R(H)=\dim _R(H)=1\ \ .\]
  \item $s=2,b=2,d=6$. We take $I:=(x_1,x_2)R\cap (x_3,x_4)R\cap
  (x_5,x_6)R$ and $H:=H^4_I(R)$. Again, because of remark
  \ref{nichtCof}, $H$ is not $I$-cofinite (it is not zero because, by Lyubeznik's result mentioned above, $\cd
  (I,R)=6-s=4$). A Mayer-Vietoris sequence argument with respect to
  the ideals $(x_1,x_2)R\cap (x_3,x_4)R$ and $(x_5,x_6)R$ shows that
  there is a canonical isomorphism $H=H^5_{(x_1,x_2,x_5,x_6)R\cap (x_3,x_4,x_5,x_6)R}(R)$. From this it is easy to see that
  \[ \injdim _R(H)=\dim _R(H)=0\ \ .\]
  \item $s=2,b=2,d=6$. We take $I:=(x_1,x_2)R\cap (x_3,x_4)R\cap
  (x_5,x_6)R$ and $H:=H^3_I(R)$. Again, because of remark
  \ref{nichtCof}, $H$ is not $I$-cofinite (we will see below that it
  is not zero). A similar Mayer-Vietoris sequence like in the
  previous example provides us with a short exact sequence
  \[ 0\to H^4_{(x_1,x_2,x_3,x_4)R}(R)\to H\to \] \[ \to H^4
  _{(x_1,x_2,x_5,x_6)R}(R)\oplus H^4_{(x_3,x_4,x_5,x_6)R}(R)\to 0\ \
  .\]
  On the other hand, the spectral sequence belonging to the composed
  functors $\Gamma _m\circ \Gamma
  _{(x_1,x_2,x_5,x_6)R}$ shows that $H^2_m(H^4
  _{(x_1,x_2,x_5,x_6)R}(R))=H^6_m(R)\neq 0$. This
  fact together with the above exact sequence
  implies $H^2_m(H)\neq 0$, therefore we have
  \[ \injdim _R(H)=\dim _R(H)=2\ \ .\]
  \item $s=2,b=3,d=7$. We take $I:=(x_1,x_2,x_3)R\cap
  (x_4,x_5,x_6)R\cap (x_7,x_1,x_2)R$ and $H:=H^5_I(R)$. Because of
  remark \ref{nichtCof}, $H$ is not $I$-cofinite. Similar arguments like
  in the first example show
  \[ \injdim _R(H)=\dim _R(H)=1\]
  (one can start e.~g. with the Mayer-Vietoris sequence belonging to
  the ideals $(x_1,x_2,x_3)R\cap (x_7,x_1,x_2)R$ and
  $(x_4,x_5,x_6)R$).
  \item $s=2,b=3,d=7$. We take $I:=(x_1,x_2,x_3)R\cap
  (x_4,x_5,x_6)R\cap (x_7,x_1,x_4)R$ and $H:=H^5_I(R)$. Because of
  remark \ref{nichtCof}, $H$ is not $I$-cofinite. The Mayer-Vietoris
  sequence with respect to the ideals $(x_1,x_2,x_3)R\cap
  (x_7,x_1,x_4)R$ and $(x_4,x_5,x_6)R$
provides us with an exact sequence
\[ H^5_{(x_1,x_2,x_3)R\cap (x_7,x_1,x_4)R}(R)\to H\to H^6_\Sigma
(R)\to H^6_{(x_1,x_2,x_3)R\cap (x_7,x_1,x_4)R}(R)\ .\] Here
\begin{eqnarray*}
\Sigma &:=&((x_1,x_2,x_3)R\cap
(x_7,x_1,x_4)R)+(x_4,x_5,x_6)R\\
&=&(x_1,\dots ,x_6)R\cap (x_1,x_4,x_5,x_6,x_7)R\ .\\
\end{eqnarray*}
Note that, by an obvious Mayer-Vietoris sequence argument, one
has
\[ H^l_{(x_1,x_2,x_3)R\cap (x_7,x_1,x_4)R}(R)=0\]
for every $l\geq 5$; therefore, one gets a canonical isomorphism
\[ H=H^6_\Sigma (R)\ .\]
The height of $\Sigma $ is five and so $H^l_\Sigma (R)=0$ for
$l\not\in \{ 5,6\} $. In particular, $H^p_{x_7R}(H^q_\Sigma (R))$ is
zero if $q\not\in \{ 5,6\} $ or if $p>1$ (i.~e. there are only four
potentially non-zero terms) and thus the convergent spectral
sequence
\[ E^{p,q}_2:=H^p_{x_7R}(H^q_\Sigma (R))\Rightarrow
H^{p+q}_{x_7R+\Sigma }(R)=H^{p+q}_{(x_1,x_4,x_5,x_6,x_7)R}(R)\]
shows $\Gamma _{x_7R}(H)=0$ (because $\Gamma
_{x_7R}(H)=E^{0,6}_2=E^{0,6}_\infty $ and even $E^6_\infty
=H^6_{(x_1,x_4,x_5,x_6,x_7)R}=0$) and $H^1_{x_7R}(H)=0$ (because
$H^1_{x_7R}(H)=E^{1,6}_2=E^{1,6}_\infty =E^7_\infty
=H^7_{(x_1,x_4,x_5,x_6,x_7)R}=0$). The two latter vanishing results
mean precisely that the canonical map $H\to H_{x_7}$ is an
isomorphism. On the other hand,
\[ (H=)H_{x_7}=H^6_\Sigma (R)_{x_7}=H^6_{(x_1,\dots ,x_6)R}(R)_{x_7}\
.\] On the other hand it is clear that
\[ \Supp _R(H^6_{(x_1,\dots ,x_6)R}(R))=\{ (x_1,\dots
,x_6)R,(x_1,\dots ,x_7)R\}\ ;\] in fact, a more precise statement is
possible:
\[ \Ext ^p_R(R/(x_1,\dots ,x_6)R,H^6_{(x_1,\dots ,x_6)R}(R))=\Ext
^{p+6}_R(R/(x_1,\dots ,x_6)R,R)\] and
\[ \Ext ^p_R(R/(x_1,\dots ,x_7)R,H^6_{(x_1,\dots ,x_6)R}(R))=\Ext
^{p+6}_R(R/(x_1,\dots ,x_7)R,R)\] ($p\in \N$ is arbitrary) are
isomorphisms that come from obvious spectral sequence arguments; the
Bass numbers of $R$ itself are known (as $R$ is Gorenstein) and so
one may conclude that a minimal injective resolution of
$H^6_{(x_1,\dots ,x_6)R}(R)$ has the form {\small \[ 0\to
H^6_{(x_1,\dots ,x_6)R}(R)\to E_R(R/(x_1,\dots ,x_6)R)\to
E_R(R/(x_1,\dots ,x_7)R)\to 0\ .\]} One may localize this exact
sequence in $x_7$ and derive canonical isomorphisms
\[ (H=)H^6_{(x_1,\dots ,x_6)R}(R)_{x_7}=E_R(R/(x_1,\dots ,x_6)R)\
.\] Now it is clear that one has
\[ \injdim _R(H)=0\neq 1=\dim _R(H)\ .\]
\end{itemize}
\end{examples}}
\begin{summary}
\label{summary} Let $R$ be a noetherian local regular ring
containing a field, $I$ an ideal of $R$ and $l$ a natural number.
Here is a summary of the previous results:

(i) If $H^l_I(R)$ is $I$-cofinite then $\injdim _R(H^l_I(R))=\dim
_R(H^l_I(R))$ holds (corollary \ref{combination}).

(ii) There are examples where $H^l_I(R)$ is not $I$-cofinite, but
$\injdim _R(H^l_I(R))=\dim _R(H^l_I(R))$ holds (e.~g. the first
example from \ref{series}).

{\sloppy (iii) There are examples where $\injdim _R(H^l_I(R))=\dim
_R(H^l_I(R))$ does not hold (e.~g. the last example from
\ref{series}, where $l=\cd (I,R)$).}
\end{summary}
\begin{question}
Do there exist a Cohen-Macaulay ideal $I$ in a regular local ring
$R$ containing a field of characteristic zero for which $\injdim
_R(H^l_I(R))=\dim _R(H^l_I(R))$ does not hold for some natural $l$?
\end{question}

{\sc Universit\"{a}t Leipzig, Fakult\"{a}t f\"{u}r Mathematik und Informatik,
Ma\-the\-ma\-tisches Institut, Augustusplatz 10/11, D-04109 Leipzig}

{\it E-mail}: hellus@math.uni-leipzig.de

\end{document}